\documentclass[11pt]{amsart}
\usepackage[colorlinks=true, pdfstartview=FitV, linkcolor=blue, citecolor=blue, urlcolor=blue, breaklinks=true]{hyperref}
\usepackage{amsmath,amsfonts,amssymb,amsthm,mathrsfs,amscd,comment,paralist,stmaryrd,etoolbox,mathtools,enumitem}
\usepackage[usenames,dvipsnames]{color}
\usepackage{booktabs}
\usepackage[all]{xy}
\usepackage{tikz}
\usetikzlibrary{decorations.pathmorphing}

\theoremstyle{plain}
\newtheorem{thm}{Theorem}[section]
\newtheorem{prop}[thm]{Proposition}
\newtheorem{lemma}[thm]{Lemma}

\theoremstyle{definition}

\newtheorem{rmk}[thm]{Remark}

\numberwithin{equation}{section}

\newcommand{\g}{\mathfrak{g}}

\newcommand{\so}{\mathfrak{so}}

\newcommand{\C}{\mathbb{C}}
\newcommand{\N}{\mathbb{N}}

\newcommand{\cint}{\int\limits}

\DeclareMathOperator{\re}{Re}

\newtoggle{comments}
\newtoggle{details}
\newtoggle{prelimnote}
\newtoggle{detailsnote}

\toggletrue{detailsnote}   

\iftoggle{comments}{%
  \newcommand{\comments}[1]{
    \begin{center}
      \parbox{6.5 in}{
        \color{red}
          {\footnotesize \textbf{Comments:} #1}
        \color{black}}
    \end{center}}
}{%
  \newcommand{\comments}[1]{}
}

\iftoggle{details}{%
  \newcommand{\details}[1]{
      \ \\
      \color{OliveGreen}
        \begin{footnotesize}
          \textbf{Details:} #1
        \end{footnotesize}
      \color{black}
      \\
  }
}{%
  \newcommand{\details}[1]{}
}

\iftoggle{prelimnote}{%
  
}{%
  
}
\allowdisplaybreaks
\begin{document}
%

\title{On the number of irreducible representations of $\so(5)$}

\author[Nayak]{Saudamini Nayak}
\address{Saudamini Nayak, Department of Mathematics, National Institute of Technology Calicut, NIT Campus P.O., Kozhikode-673 601, India}
\email{saudamini@nitc.ac.in}

\author[Meher]{N. K. Meher}
\address{Nabin Kumar Meher, National Institute of Technology Raipur, G.E. Road, Raipur, Chhattisgarh, 492010, India. }
\email{mehernabin@gmail.com, nkmeher.maths@nitrr.ac.in}

\author[Rout]{S. S. Rout}
\address{Sudhansu Sekhar Rout, Department of Mathematics, National Institute of Technology Calicut, NIT Campus P.O., Kozhikode-673 601, India}
\email{lbs.sudhansu@gmail.com; sudhansu@nitc.ac.in}

\thanks{}

\begin{abstract}
Let $d(n)$ be the divisor function and it is well known that  $\sum_{1\leq n \leq x}d(n) = x\log x+(2\gamma-1)x +\mathcal{O}\left(x^{\theta+\epsilon}\right)$
where $\gamma$ is the Euler constant, $\epsilon>0$ and $1/4<\theta<1/3$.  In this paper, we obtain an asymptotic formula  for the number of irreducible representations of $\mathfrak{so}(5)$. More precisely,
the irreducible representations of the Lie algebra $\mathfrak{so}(5)$ are a family of representations of dimension $jk(j+k)(j+2k)/6$ for $j, k\in \mathbb{N}_{0}$ and suppose that $\varrho_{\mathfrak{so}(5)}(n)$ is the number of irreducible $\mathfrak{so}(5)$ representations of dimension $n$. We obtain an asymptotic formula for the summatory function $\sum_{1\leq n \leq x}\varrho_{\mathfrak{so}(5)}(n)$. 
\end{abstract}

\subjclass[2020]{17B15, 11N45, 11N56}
\keywords{Witten zeta function, asymptotic formula, representations of Lie algebra, }

\maketitle
\thispagestyle{empty}

\setcounter{tocdepth}{1}


%
\section{Introduction}
It is well-known that suitably defined zeta and $L$-functions and their special values often play significant roles in many areas of mathematics. Witten \cite{witten} studied one variable zeta functions attached to various Lie algebras and related their special values to the volumes of certain moduli spaces of vector bundles of curves. The Witten zeta-function associated with a complex semisimple Lie algebra $\g$ is defined as 
\begin{equation}\label{eq1}
    \zeta_{\g}(s) = \sum_{\rho} \frac{1}{(\dim \rho)^s},
\end{equation}
where the summation runs over all finite dimensional irreducible representations $\rho$ of $\g$. The above definition is due to Zagier \cite{zagier}, and the values $\zeta_{\g}(2k)$  for positive integers $k$ were first studied by Witten \cite{witten}. In \cite{zagier} Zagier proved that
\begin{equation}\label{eq2}
    \zeta_{\g}(2k) \in \mathbb{Q} \pi^{kl}\quad (k\in \N),
\end{equation}
where $l$ is the number of positive roots of $\g$, using Witten’s result \cite{witten}. Some explicit formulas for $\zeta_{\g}(2k)\, (k \in \N)$ were given by Mordell \cite{mordell}, Zagier \cite{zagier}, Subbarao and Sitaramachandrarao \cite{subbarao-etal} and Gunnells and Sczech \cite{gunnels-sczech}. More recently Matsumoto and Tsumura \cite{matsumoto-semisimple1} and Komori et al. \cite{matsumoto-semisimple2} introduced the multi-variable Witten zeta-functions associated with semisimple Lie algebras, and studied their analytical and arithmetical properties (see also \cite{matsumoto-semisimple3, matsumoto-semisimple4}).

The Lie algebra $\mathfrak{su}(2)$ has (up to equivalence) one irreducible representation $V_k$ of each dimension $k\in \mathbb{N}$. Each $n$-dimensional representation $\oplus_{k=1}^{\infty}r_kV_k$ corresponds to a unique partition
\begin{equation}\label{eq1.3}
    n = \lambda_1+\lambda_2+\cdots+\lambda_r, \quad 1\leq \lambda_1\leq \lambda_2\leq \cdots\lambda_r\leq n
\end{equation}
such that $r_k$ denotes the number of parts in the partition that equal to $n$. So, the number of $n$-dimensional representations equals $p(n)$, the number of integer partitions of $n$. The partition function has no elementary closed formula, nor does it satisfy any finite order recurrence. However, with $p(0)=1$, its generating function has the following classical product expansion:
\begin{equation}\label{eq1.4}
    \sum_{n=0}^{\infty}p(n)q^n = \prod_{n=1}^{\infty}\frac{1}{1-q^n}.
\end{equation}
In \cite{hardy-ramanujan}, Hardy and Ramanujan initiated the analytic study of $p(n)$ and used circle method to find the asymptotic formula for the partition function 
\begin{equation}
p(n) \sim \frac{1}{4\sqrt{3}n}e^{\pi \sqrt{2n/3}} \quad \mbox{as}\;\; n\to \infty.
\end{equation}
The next case is the Lie algebra $\mathfrak{su}(3)$, whose irreducible representations $W_{j,k}$ indexed by pairs of positive integers (see \cite[Theorem 6.27]{hall}. Like in the case of $\mathfrak{su}(2)$, a general $n$-dimensional representation decomposes into a sum of these $W_{j,k}$, again each with some multiplicity. So analogous to \eqref{eq1.4}, the numbers $r_{\mathfrak{su}(3)}(n)$ of $n$-dimensional representations, have the generating function
\begin{equation}\label{eq1.5}
    \sum_{n\geq 0}r_{\mathfrak{su}(3)}(n)q^n = \prod_{j, k\geq 1}\frac{1}{1-q^{\frac{jk(j+k)}{2}}}
\end{equation}
with $r_{\mathfrak{su}(3)}(0) =1$. Romik in \cite{romik2017} proved that as $n\to \infty$ 
\[r_{\mathfrak{su}(3)}(n)\sim \frac{K}{n^{\frac{2}{3}}} \exp\left(A_1n^{\frac{2}{5}} - A_2n^{\frac{3}{10}}-A_3n^{\frac{1}{5}}- A_4n^{\frac{1}{10}}\right).\]
Recently, Bringmann and Franke \cite{Bringmann2023} improved this as follows:
\begin{align*}
 r_{\mathfrak{su}(3)}(n)&\sim \frac{C_0}{n^{\frac{2}{3}}} \exp\left(A_1n^{\frac{2}{5}} + A_2n^{\frac{3}{10}}+A_3n^{\frac{1}{5}} + A_4n^{\frac{1}{10}}\right)\left(1+\sum_{j=1}^N\frac{C_j}{n^{\frac{j}{10}}}+ \mathcal{O}_N\left(n^{-\frac{N+1}{10}}\right)\right)
\end{align*}
when $n\to \infty$.
This framework also generalizes to the Witten zeta function for $\mathfrak{so}(5)$. It is known that the finite-dimensional representations of $\mathfrak{so}(5)$ can be doubly indexed as $(T_{j,k})_{j,k \in \mathbb{N}}$ with $\dim T_{j,k} = \frac{1}{6}jk(j + k)(j + 2k)$. A general $n$-dimensional representation decomposes as a sum of these $T_{j,k}$, each with some multiplicity. Therefore, one can have
\begin{equation}\label{eq1.5}
    \sum_{n\geq 0}r_{\mathfrak{so}(5)}(n)q^n = \prod_{j, k\geq 1}\frac{1}{1-q^{\frac{jk(j+k)(j+2k)}{6}}}.
\end{equation}
In \cite{Bridges2024}, Bridges et al., derived an assymptotic formula for $r_{\mathfrak{so}(5)}(n)$ for any $N\in \mathbb{N}$ and as $n\to \infty$
\begin{align*}
  r_{\mathfrak{so}(5)}(n)&\sim  \frac{C}{n^{\frac{7}{12}}} \exp\left(A_1n^{\frac{1}{3}} + A_2n^{\frac{2}{9}}+A_3n^{\frac{1}{9}} + A_4\right) \left(1+\sum_{j=2}^{N+1}\frac{B_j}{n^{\frac{j-1}{9}}}+ \mathcal{O}_N\left(n^{-\frac{N+1}{9}}\right)\right).
\end{align*}
Let $\varrho_{\mathfrak{su}(3)}(n)$ denote the number of irreducible $\mathfrak{su}(3)$ representations of dimension $n$, that is,
\begin{equation}
  \varrho_{\mathfrak{su}(3)}(n) = \sum_{\substack{j, k\geq 1 \\ \frac{jk(j+k)}{2} =n}}  1.
\end{equation}
Note that this is a highly oscillatory function and is often $0$ and it is similar to the divisor function \[d(n) = \sum_{\substack{j, k\geq 1 \\ jk=n}}1.\] 
We know that Dirichlet hyperbola method gives 
\begin{equation}
    \sum_{1\leq n \leq x}d(n) = x\log x+(2\gamma-1)x +\mathcal{O}(\sqrt{x}),
\end{equation}
where $\gamma$ is the Euler-Mascheroni constant. Using a variant of the the hyperbola method, Bridges et al. \cite{Bridges-Bringmann2024}, obtained an assymptotic expansion for the function $\varrho_{\mathfrak{su}(3)}(n)$ as follows:
\begin{equation}
   \sum_{1\leq n \leq x} \varrho_{\mathfrak{su}(3)}(n) = \frac{2^{\frac{2}{3}\sqrt{3}\Gamma(\frac{1}{3})^3}}{4\pi}x^{\frac{2}{3}}+ 2^{\frac{3}{2}}\zeta\left(\frac{1}{2}\right)\sqrt{x}+ \mathcal{O}(x^{\frac{1}{3}}),
\end{equation}
as $x\to \infty$.
In  the same paper \cite{Bridges-Bringmann2024}, they pose the question of establishing an asymptotic formula for $\displaystyle \sum_{\substack{j, k\geq 1 \\ p(j, k)\leq x}} 1$, where $p(x, y)$ is a homogeneous polynomial in $\mathbb{Q}[x, y]$ taking integer values. In this paper, we deduce an asymptotic expansion for the number of irreducible $\mathfrak{so}(5)$-representaion of dimension $n$, that is,
\begin{equation}
\sum_{1\leq n\leq x}  \varrho_{\mathfrak{so}(5)}(n) =\sum_{1\leq n\leq x} \sum_{\substack{j, k\geq 1 \\ \frac{jk(j+k)(j+2k)}{6} =n}}  1.
\end{equation}
Precisely we prove the following.
\begin{thm}\label{thm1}
We have, as $x\to \infty$
  \[ \sum_{1\leq n\leq x}\varrho_{\mathfrak{so}(5)}(n) = \frac{\sqrt{3}\Gamma(\frac{1}{4})^2}{4\sqrt{\pi}} x^{\frac{1}{2}}+ K\cdot x^{\frac{1}{3}} + \mathcal{O}(x^{\frac{1}{4}}),\]
  where
  \[K:=\left(\left(\zeta\left(\frac{1}{3} \right)+\frac{1}{2}\right) (10 \cdot 6^{1/3}+2\cdot 3^{1/3}) - 3^{4/3} (2^{-2/3}+2^{-1}) + 3^{\frac{1}{3}}(2+2^{\frac{4}{3}})\right).\]
\end{thm}
    
The organizations of the paper is as follows. In Section \ref{prelim}, we describe the multi-variable Witten zeta-functions associated with semisimple Lie algebra $\so(5)$. Next we discuss the analytic continuation of the series $\zeta_{\so(5)}(s)$ and the location of its poles and residues. Finally, in Section \ref{secproof}, we proof the main theorem.

\section{Preliminaries}\label{prelim}
\subsection{Witten zeta function associated to $\so(5)$}
In this section, we describe the multi-variable Witten zeta-functions associated with semisimple Lie algebra $\so(5)$ in terms of roots and weights for the corresponding root systems using Weyl’s dimension formula \cite{knapp}. We follow Matsumoto and Tsumura \cite{matsumoto-semisimple1}  and Komori et al. \cite{matsumoto-semisimple2} for the following definitions. 

The fundamental system for $\so(5)$ is given by $\Pi:= \{\alpha_1=\epsilon_1-\epsilon_2, \alpha_2= \epsilon_2\}$ and the list of positive roots are:
\begin{equation}\label{eq3}
    \Delta_{+}:= \{\epsilon_1= \alpha_1+\alpha_2, \epsilon_2 = \alpha_2, \epsilon_1-\epsilon_2= \alpha_1, \epsilon_1+\epsilon_2= \alpha_1+2\alpha_2\}.
\end{equation}
The fundamental coroots are given by $\{\alpha_1^{\vee}=e_1-e_2, \alpha_2^{\vee} = 2e_2\}$ and hence the positive coroots are $\{2e_1, 2e_2, e_1-e_2, e_1+e_2\}$. Let $\lambda_1, \lambda_2$ be two fundamental weights satisfying $<\alpha_i^{\vee}, \lambda_j> = \lambda_j(\alpha_i) = \delta_{ij}$ where $\delta_{ij}$ is Kronecker's delta function and $i, j=1, 2$. Any dominant weight can be written as \[\lambda = n_1\lambda_1+n_2\lambda_2 \quad n_1, n_2\in \mathbb{N}\cup\{0\}.\] The lowest strongly dominant is of the form $\delta = \lambda_1+\lambda_2$. Let $d_{\lambda}$ be the dimension of the representation space corresponding to the weight $\lambda$. By Weyl's dimension formula, we have the following.
\begin{align*}
    d_{\lambda} &= \prod_{\alpha\in \Delta_{+}} \frac{<\alpha^{\vee}, \lambda+\delta>}{<\alpha^{\vee}, \delta>}\\
    & = \prod_{\alpha\in \Delta_{+}} \frac{<\alpha^{\vee}, (n_1+1)\lambda_1+(n_2+1)\lambda_2>}{<\alpha^{\vee}, \lambda_1+\lambda_2>} \\
    & = (n_1+1) (n_2+1) (n_1+n_2+2) (n_1+2n_2+3)/6.
\end{align*}
Now set $n_1+1=m, n_2+1=n$. Then from \eqref{eq1}, Witten zeta function is given by
\begin{equation}\label{eq4}
  \zeta_{\so(5)}(s) = 6^s\sum_{n=1}^{\infty}\sum_{m=1}^{\infty} \frac{1}{m^s n^s (m+n)^s (m+2n)^s}.  
\end{equation}
To obtain the precise region of convergence of the Dirichlet series given in \eqref{eq4}, we rewrite \eqref{eq4} as
\begin{align}
\begin{split}\label{eq5}
    \zeta_{\so(5)}(s) &= 6^s \left(\sum_{n>m\geq 1}^{\infty} + \sum_{m>n\geq 1}^{\infty}+ \sum_{n = m\geq 1}^{\infty}\right)\sum_{n=1}^{\infty}\sum_{m=1}^{\infty} \frac{1}{m^s n^s (m+n)^s (m+2n)^s}\\
    &= 2\cdot 6^s\sum_{n>m\geq 1}^{\infty}  \frac{1}{m^s n^s (m+n)^s (m+2n)^s}+ \zeta(4s)\\
    & = 2\cdot 6^s\sum_{n>m\geq 1}^{\infty}  \frac{1}{m^{s} n^{3s} (1+m/n)^s (2+m/n)^s}+ \zeta(4s).
\end{split}
\end{align}
Note that the series $\zeta(4s)$ converges if and only if $\re(s)>1/4$ and since 
\[6^{-s}< \sum_{n>m\geq 1}^{\infty}\frac{1}{(1+m/n)^s (2+m/n)^s} <2^{-s},\] 
then the first sum in \eqref{eq5} converges absolutely if and only if the series 
\begin{equation}\label{eq6}
    \sum_{n>m\geq 1}^{\infty}\frac{1}{|n^{3s} m^s|} = \sum_{n>m\geq 1}^{\infty}\frac{1}{n^{3\sigma} m^{\sigma}}, \quad \mbox{where}\;\; s = \sigma+it
\end{equation}
 converges absolutely. In comparision with the integral 
 \[\int_{1}^{\infty}\int_{x}^{\infty} \frac{1}{x^{3\sigma}y^{\sigma}}dx dy = \frac{1}{\sigma-1} \int_{1}^{\infty} x^{1-4\sigma}dx,\]
we conclude that the series in \eqref{eq6} converges absolutely and uniformly on compacts if and only if  $\re(s)>1/2$. The above discussion proves the following:
\begin{prop}
The series \eqref{eq4} converges absolutely in the domain $\{s\in \mathbb{C}\mid \re(s) >1/2\}$ and defines a holomorphic function in that region.
\end{prop}

\subsection{Analytic continuation of $\zeta_{\so(5)}(s)$}
Here we will see the analytic continuation of the series $\zeta_{\so(5)}(s)$ to a meromorphic function on $\mathbb{C}$. This fact can be obtained from \cite{matsumoto-bonner}. To obtain analytic continuation of $\zeta_{\so(5)}(s)$, we need the {\em Mordell–Tornheim} zeta function $\zeta_{MT}$. Let $\re(s_j)>1, (j=1, 2, 3)$ and define
\begin{equation}
    \zeta_{MT}(s_1, s_2, s_3) = \sum_{n=1}^{\infty}\sum_{m=1}^{\infty} \frac{1}{m^{s_1} n^{s_2} (m+n)^{s_3}}.
\end{equation}
\begin{lemma}[\cite{matsumoto}]
The function $\zeta_{MT}(s_1, s_2, s_3)$ can be meromorphically continued to the whole $\mathbb{C}^3$-space, and its singularities are only on the subsets of $\mathbb{C}^3$ defined by one of the equations 
\begin{equation}
    s_1+s_3 =1-l, s_2+s_3 =1-l,\,(l\in \mathbb{N}_0) \;\; \mbox{or}\;\; s_1+s_2+s_3 =2.
\end{equation}
\end{lemma}
We recall Mellin-Barnes integrals \cite[ p.~91]{paris-kaminski}, namely
\begin{equation}
\label{eq:mellin-barnes}
\Gamma(s) (1+\lambda)^{-s} = \frac{1}{2\pi i} \cint_{\alpha-i\infty}^{\alpha+i\infty} \Gamma(s+z)\Gamma(-z) \lambda^z\,dz,
\end{equation}
where  $\lambda \in \C\setminus (-\infty,0]$ and $-\re(s)<\alpha<0$.  Setting $\lambda=n/(m+n)$, multiplying by $m^{-s}n^{-s}(m+n)^{-2s}$ and summing over $n,m\in \mathbb{N}$ gives
\begin{align}
\Gamma(s) \zeta_{\so(5)}(s) &= 
\Gamma(s) 6^s \sum_{n,m=1}^\infty m^{-s} n^{-s}\left(m+n\right)^{-2s}\left(1+\frac{n}{m+n}\right)^{-s}
\nonumber \\&= \frac{6^s}{2\pi i} \sum_{n,m=1}^\infty m^{-s} n^{-s} (m+n)^{-2s}\cint_{\alpha-i\infty}^{\alpha+i\infty} \Gamma(s+z)\Gamma(-z) (m+n)^{-z} n^{z} \,dz
\nonumber \\ &=\frac{6^s}{2\pi i}\cint_{\alpha-i\infty}^{\alpha+i\infty}\Gamma(s+z)\Gamma(-z) \ \zeta_{MT}(s, s-z, 2s+z)  \,dz.   
\label{eq:wzeta-integral-rep1}
\end{align}
where $-\re(s)<\alpha<0$.  Let $L$ be a large positive integer, and put $\Phi(s) = (4s-2)\prod_{\ell =0}^L (3s-1+\ell)$. Then 
\begin{equation}\label{}
\zeta_{\so(5)}(s) = \Phi(s)^{-1} I,
\end{equation} 
where 
\begin{equation}\label{eq2.8}
I=  \frac{1}{2\pi i}\cint_{\alpha-i\infty}^{\alpha+i\infty}\frac{\Gamma(s+z)\Gamma(-z)}{\Gamma(s)}  \zeta_{MT}(s, s-z, 2s+z)\Phi(s)  \,dz.
\end{equation}
We shift the path of integration to $\re(z) = M -\epsilon$, where $M$ is a large positive integer and $\epsilon$ is a small positive number. Counting the residues at $z = k\ (0 \leq k\leq  M -1)$, we obtain
\begin{align}\label{eq2.9}
\begin{split}
   \zeta_{\so(5)}(s) &= \frac{6^s}{\Gamma(s)}\sum_{k=0}^{M-1}\frac{(-1)^{k+1}}{k!}\frac{\Gamma(s+k)}{\Gamma(s)}\zeta_{MT}(s, s-k, 2s+k)+\Phi(s)^{-1}I' \\
   &= \frac{6^s}{\Gamma(s)}\sum_{k=0}^{M-1}\binom{-s}{k}\zeta_{MT}(s, s-k, 2s+k)+\Phi(s)^{-1}I' 
\end{split}
\end{align}
where
\begin{equation}\label{eq2.10}
  I'= \frac{1}{2\pi i}\cint_{M-\epsilon-i\infty}^{M-\epsilon+i\infty} \frac{\Gamma(s+z)\Gamma(-z)}{\Gamma(s)}  \zeta_{MT}(s, s-z, 2s+z)\Phi(s)\,dz. 
\end{equation}
Hence, the integral $I'$ can be continued holomorphically to the region
\[\mathcal{D}_{M,L}:=\{s\in \mathbb{C}\mid \re(s) >-M+\epsilon, \re(3s)>-L, \;\re(3s)> 1-M+\epsilon\}\]
because in this region the poles of the integrand are not on the path of integration. Therefore, \eqref{eq2.9} gives the meromorphic continuation of $\zeta_{\so(5)}(s)$ to $\mathcal{D}_{M,L}$, and the candidates for singularities in this region are
\begin{equation}\label{eq2.11}
   3s = 1-\ell,\; 4s = 2 , \quad \ell \in \mathbb{N}_0. 
\end{equation}
\begin{prop}[\cite{Bridges2024}]\label{prop-analcont}
The function $\zeta_{\so(5)}(s)$ can be analytically continued to a meromorphic function on $\mathbb{C}$ whose positive poles are simple and occur for $s\in \left\{\frac{1}{2}, \frac{1}{3}, -\frac{1}{3}, -\frac{2}{3}, \ldots\right\}$.
\end{prop}
Furthermore, we also need certain residues of $\zeta_{\so(5)}(s)$. For a detailed proof of the following result, refer to \cite[Proposition 5.16]{Bridges2024}.
\begin{prop}\label{prop-poles}
    The poles of $\zeta_{\so(5)}(s)$ are precisely $\{\frac{1}{2}\}\cup \{\frac{d}{3}\not \in \mathbb{Z}\mid d\leq 1, d\;\;\mbox{odd}\}$. We have 
    \[\mbox{Res}_{s=\frac{1}{2}}\zeta_{\so(5)}(s) = \frac{\sqrt{3}\Gamma(\frac{1}{4})^2}{8\sqrt{\pi}}.\] Moreover, for $d\in \mathbb{Z}_{\leq 1} \setminus \{-3\mathbb{N}_{0}\}$
    \[\mbox{Res}_{s=\frac{d}{3}}\zeta_{\so(5)}(s) = \frac{3^{\frac{d}{3}-\frac{3}{2}}\pi\Gamma(\frac{d}{6})\zeta(\frac{4d}{3}-1)}{2^{\frac{d}{3}-1}(1-d)!\Gamma(\frac{d}{3})^2\Gamma(\frac{d}{2})}\left(\frac{d}{3}\right)(1+2^{\frac{2d}{3}-1}).\] 
\end{prop}

\section{Proof of Theorem \ref{thm1}}\label{secproof}
At first, we note that the asymptotic main term of $\sum_{1\leq n\leq x}\varrho_{\mathfrak{so}(5)}(n)$ can be derived by the analytic properties of $\zeta_{\so(5)}(s)$ along with a standard Tauberian theorem. Observe that we rewrite
\[\zeta_{\so(5)}(s)=\sum_{j,k \geq 1}\frac{6^s}{j^sk^s(j+k)^s(j+2k)^{s}} = \sum_{n\geq 1} \frac{\varrho_{\mathfrak{so}(5)}(n)}{n^s}.\]
By Proposition \ref{prop-analcont}, $\zeta_{\so(5)}(s)$ has a meromorphic continuation to $\mathbb{C}\setminus \{\frac{1}{2}\}\cup \{\frac{d}{3}\not \in \mathbb{Z}\mid d\leq 1, d\;\;\mbox{odd}\}$. Then from the Ikehera-Wiener Tauberian theorem (see \cite[Chapter 3]{murty}), we get
\begin{equation}\label{eq3.1}
   \sum_{1\leq n\leq x}\varrho_{\mathfrak{so}(5)}(n) \sim \frac{\mbox{Res}_{s=\frac{1}{2}}\zeta_{\so(5)}(s)}{\frac{1}{2}}x^{\frac{1}{2}}=\frac{\sqrt{3}\Gamma(\frac{1}{4})^2}{4\sqrt{\pi}} x^{\frac{1}{2}}, \quad \mbox{as}\;\; x\to \infty. 
\end{equation}
We write
\begin{equation}\label{mainsum}
  \sum_{1\leq n\leq x}\varrho_{\mathfrak{so}(5)}(n)= \sum_{1\leq N\leq x} \sum_{\substack{m, n\geq 1 \\ \frac{mn(m+n)(m+2n)}{6} =N}}  1= \sum_{\substack{m, n\geq 1 \\ mn(m+n)(m+2n) \leq 6x\\ mn(m+n)(m+2n)\equiv 0\pmod{6}}}  1.  
\end{equation}
to obtain the next term in the asymptotic expansion. Since $mn(m+n)(m+2n)\equiv 0\pmod{2}$ and $mn(m+n)(m+2n)\equiv 0\pmod{3}$, we can conclude that $mn(m+n)(m+2n)\equiv 0\pmod{6}$ is always satisfied. Observe that the sum in \eqref{mainsum} counts lattice points in the $(m, n)$-plane between $m=1, n=1$ and the curves given by $m=T_1$ and $n = T_2$ where
\begin{equation}\label{curve1}
    T_1 := n\left[\left(\frac{3x}{n^4}+\frac{1}{2}\sqrt{\left(\frac{6x}{n^4}\right)^2-\frac{4}{27}}\right)^{1/3} + \left(\frac{3x}{n^4}-\frac{1}{2}\sqrt{\left(\frac{6x}{n^4}\right)^2-\frac{4}{27}}\right)^{1/3}-1\right]
\end{equation}and 
\begin{equation}\label{curve2}
    T_2 := m\left[\left(\frac{3x}{2m^4}+\frac{1}{2}\sqrt{\left(\frac{3x}{m^4}\right)^2-\frac{1}{432}}\right)^{1/3} + \left(\frac{3x}{2m^4}-\frac{1}{2}\sqrt{\left(\frac{3x}{m^4}\right)^2-\frac{1}{432}}\right)^{1/3}-\frac{1}{2}\right].
\end{equation}
Note that $m$ and and $n$ given in \eqref{curve1} and \eqref{curve2} are positive solutions to the cubic equation \[mn(m+n)(m+2n) =6x.\]  Then, using the hyperbola method and the notation $\{t\}:= t-\lfloor{t\rfloor}$, we rewrite \eqref{mainsum} as follows.
\begin{align}\label{align3.2}
\begin{split}
  \sum_{1\leq n\leq x}\varrho_{\mathfrak{so}(5)}(n) &= \sum_{1\leq n\leq x^{\frac{1}{4}}} \sum_{1\leq m\leq T_1} 1 + \sum_{1\leq m\leq x^{\frac{1}{4}}} \sum_{1\leq n\leq T_2} 1 - \lfloor{x^{\frac{1}{4}}\rfloor}^2\\
  &= \sum_{1\leq n\leq x^{\frac{1}{4}}}\lfloor{T_1\rfloor} + \sum_{1\leq m\leq x^{\frac{1}{4}}}\lfloor{T_2\rfloor}- \lfloor{x^{\frac{1}{4}}\rfloor}^2\\
  &= \sum_{1\leq n\leq x^{\frac{1}{4}}}\lfloor{T_1\rfloor} + \sum_{1\leq m\leq x^{\frac{1}{4}}}\lfloor{T_2\rfloor}-x^{\frac{1}{2}}+ \mathcal{O}(x^{\frac{1}{4}}).
  \end{split}
\end{align}
First consider
\begin{align}\label{align3.6}
\begin{split}
   \sum_{1\leq n\leq x^{\frac{1}{4}}}\lfloor{T_1\rfloor} &= \sum_{1\leq n\leq x^{\frac{1}{4}}}n\left(\left(\frac{3x}{n^4}+\frac{1}{2}\sqrt{\left(\frac{6x}{n^4}\right)^2-\frac{4}{27}}\right)^{1/3} + \left(\frac{3x}{n^4}-\frac{1}{2}\sqrt{\left(\frac{6x}{n^4}\right)^2-\frac{4}{27}}\right)^{1/3}\right)\\
   &-\sum_{1\leq n \leq x^{\frac{1}{4}}}n + \mathcal{O}(x^{\frac{1}{4}})=\sum_{1\leq n\leq x^{\frac{1}{4}}} f_1(n)- \frac{ \lfloor{ x^{\frac{1}{4}} \rfloor}  (\lfloor{ x^{\frac{1}{4}} \rfloor}+1) }{2} +\mathcal{O}(x^{\frac{1}{4}})\\
   &= \sum_{1\leq n\leq x^{\frac{1}{4}}} f_1(n) - \frac{1}{2} x^{\frac{1}{2}} + \mathcal{O}(x^{\frac{1}{4}}),
   \end{split}
\end{align} where 
\begin{align*}
	f_1(n) &= n\left[\left(\frac{3x}{n^4}+\frac{1}{2}\sqrt{\left(\frac{6x}{n^4}\right)^2-\frac{4}{27}}\right)^{1/3} + \left(\frac{3x}{n^4}-\frac{1}{2}\sqrt{\left(\frac{6x}{n^4}\right)^2-\frac{4}{27}}\right)^{1/3}\right].
\end{align*}
We apply the Abel partial summation formula \cite[p.4]{tenenbaum} to approximate the sum $\sum_{1\leq n\leq x^{\frac{1}{4}}} f_1(n)$. To do this, we set $a_n=1, A(t) = \sum_{1\leq n\leq t}a_n$ and 
\begin{align*}
f_1(t) &= t\left[\left(\frac{3x}{t^4}+\frac{1}{2}\sqrt{\left(\frac{6x}{t^4}\right)^2-\frac{4}{27}}\right)^{1/3} + \left(\frac{3x}{t^4}-\frac{1}{2}\sqrt{\left(\frac{6x}{t^4}\right)^2-\frac{4}{27}}\right)^{1/3}\right].
\end{align*}
Thus, 
\begin{align}\label{eq3.6a}
    \begin{split}
 \sum_{1\leq n\leq x^{\frac{1}{4}}} f_1(n)&= A(x^{1/4})f_1(x^{1/4})-\int_{1}^{x^{1/4}} A(t)f_1'(t)dt =-\int_{1}^{x^{1/4}} \lfloor{t\rfloor}f_1'(t)dt+ \mathcal{O}(x^{\frac{1}{4}})\\
 &= -\int_{1}^{x^{1/4}} tf_1'(t)dt + \int_{1}^{x^{1/4}} \{t\}f_1'(t)dt+\mathcal{O}(x^{\frac{1}{4}})\\
 & = -{x^{1/4}}f_1(x^{1/4})+f_1(1) + \int_{1}^{x^{1/4}}f_1(t)dt+\int_{1}^{x^{1/4}} \{t\}f_1'(t)dt+\mathcal{O}(x^{\frac{1}{4}})\\
 &= f_1(1)+\int_{1}^{x^{1/4}}f_1(t)dt+ \int_{1}^{x^{1/4}} \{t\}f_1'(t)dt+\mathcal{O}(x^{\frac{1}{4}}). 
 \end{split}
\end{align}
Setting $\sqrt{r}:=  \left(\frac{3x}{t^4}\right)^{-1}\sqrt{\left(\frac{3x}{t^4}\right)^2-\frac{1}{27}}$, we get
\begin{align}\label{eqnewf1}
\begin{split}
\frac{f_1(t)}{t} &=\left(\frac{3x}{t^4}\right)^{1/3}[(1+ \sqrt{r})^{1/3}+(1-\sqrt{r})^{1/3}]\\
& = 2\left(\frac{3x}{t^4}\right)^{1/3} \sum_{m=0}^{\infty} \binom{1/3}{2m} r^{m}\\
&= 2 \left(\frac{3x}{t^4}\right)^{1/3}\left[1-\frac{1}{9}r-\frac{10}{243}r^2-\cdots\right] =  2^{4/3} \left(\frac{3x}{t^4}\right)^{1/3} + \mathcal{O} \left( \left(3x\right)^{-5/3} t^{20/3}\right).
\end{split}
\end{align}
Putting $t=1$ in \eqref{eqnewf1} and  then using it in \eqref{eq3.6a}, we infer
\begin{equation}\label{eq3.6}
\sum_{1\leq n\leq x^{\frac{1}{4}}} f_1(n)= 2^{4/3}(3x)^{1/3}+\int_{1}^{x^{1/4}}f_1(t)dt+ \int_{1}^{x^{1/4}} \{t\}f_1'(t)dt+\mathcal{O}(x^{\frac{1}{4}}). 
\end{equation}
Note that we use  \eqref{eqnewf1} to evaluate $f_1(1)$ in the last equality. Upon changing the variables $t\to (3x)^{\frac{1}{4}}t$, the first integral $\int_{1}^{x^{1/4}}f(t)dt$ is
\begin{align}
\begin{split}
    \int_{1}^{x^{\frac{1}{4}}}f_1(t)dt &= \int_{1}^{x^{\frac{1}{4}}}t\left[\left(\frac{3x}{t^4}+\frac{1}{2}\sqrt{\left(\frac{6x}{t^4}\right)^2-\frac{4}{27}}\right)^{1/3} + \left(\frac{3x}{t^4}-\frac{1}{2}\sqrt{\left(\frac{6x}{t^4}\right)^2-\frac{4}{27}}\right)^{1/3}\right]dt\\
    &= (3x)^{\frac{1}{2}}\int_{\frac{1}{(3x)^{1/4}}}^{\frac{1}{3^{1/4}}}t\left[\left(\frac{1}{t^4}+\sqrt{\left(\frac{1}{t^4}\right)^2-\frac{1}{27}}\right)^{1/3} + \left(\frac{1}{t^4}-\sqrt{\left(\frac{1}{t^4}\right)^2-\frac{1}{27}}\right)^{1/3}\right]dt\\
    &=:(3x)^{\frac{1}{2}}G_1\left(\frac{1}{(3x)^{\frac{1}{4}}}\right)
\end{split}
\end{align}
where \[G_1(y):= \int_{y}^{\frac{1}{3^{1/4}}}t\left[\left(\frac{1}{t^4}+\sqrt{\left(\frac{1}{t^4}\right)^2-\frac{1}{27}}\right)^{1/3} + \left(\frac{1}{t^4}-\sqrt{\left(\frac{1}{t^4}\right)^2-\frac{1}{27}}\right)^{1/3}\right]dt.\]
Rewrite 
\[G_1(y)= \int_{y}^{\frac{1}{3^{1/4}}} t^{-1/3} \left[\left(1+\sqrt{1-\frac{t^8}{27}}\right)^{1/3} + \left(1-\sqrt{1-\frac{t^8}{27}}\right)^{1/3}\right]dt\]
and then simplifying, we get the following.
\begin{align*}
G_1(0)& = \int_{0}^{y} t^{-1/3} \left[\left(1+\sqrt{1-\frac{t^8}{27}}\right)^{1/3} + \left(1-\sqrt{1-\frac{t^8}{27}}\right)^{1/3}\right]dt + G_1(y)\\
&= G_1(y)+ \int_{0}^{y} t^{-1/3} H\left(1-\frac{t^8}{27}\right)dt
\end{align*}
where 
\begin{equation}\label{eqnewhz}
H(z) := (1+\sqrt{z})^{1/3}+ (1 - \sqrt{z})^{1/3}.
\end{equation} 
Using the binomial expansion, we deduce that
\begin{equation}\label{binom1}
H(z) = 2-\frac{2}{9}z-\frac{20}{243}z^2-\cdots
\end{equation}
and this implies
\[G_1(y) = G_1(0) - 3\cdot 2^{-\frac{2}{3}} y^{\frac{2}{3}} +\mathcal{O}(y^{\frac{26}{3}}).\]
Hence,
\begin{align}
\begin{split}\label{eq3.8}
    \int_{1}^{x^{\frac{1}{4}}}f_1(t)dt &=(3x)^{\frac{1}{2}}\left[G_1(0) - 3\cdot 2^{-\frac{2}{3}} \left(\frac{1}{(3x)^{\frac{1}{4}}}\right)^{\frac{2}{3}} +\mathcal{O}\left(\left(\frac{1}{(3x)^{\frac{1}{4}}}\right)^{\frac{26}{3}}\right)\right]\\
    &= (3x)^{\frac{1}{2}}G_1(0) - 3^{\frac{4}{3}}\cdot 2^{-\frac{2}{3}} x^{\frac{1}{3}} + \mathcal{O}(x^{-\frac{5}{3}}).
    \end{split}
\end{align}
To compute $\int_{1}^{x^{1/4}} \{t\}f_1'(t)dt$, we write $f_1(t) = t(g_1(t)+g_2(t))$ and so
 \[f_1'(t) = (g_1(t)+g_2(t)) + t (g_1'(t)+g_2'(t)) = \frac{f_1(t)}{t} + t (g_1'(t)+g_2'(t)). \]
Then
\begin{equation}\label{eq3.9}
\int_{1}^{x^{1/4}} \{t\}f_1'(t)dt = \int_{1}^{x^{1/4}} \{t\}\frac{f_1(t)}{t}dt+ \int_{1}^{x^{1/4}} t\{t\}(g_1'(t)+ g_2'(t))dt.
\end{equation}
Also, from \eqref{eqnewf1}, we get
\begin{align*}
    \int_{1}^{x^{1/4}} \{t\}\frac{f_1(t)}{t}dt& = 2^{4/3}(3x)^{1/3}\int_{1}^{x^{1/4}} \frac{\{t\}}{t^\frac{4}{3}}dt + \mathcal{O} \left( \left(3x\right)^{-5/3} \int_{1}^{x^{1/4}} \{t\} t^{20/3} dt \right) \\
    & = 2^{4/3}(3x)^{1/3}\left[\int_{1}^{\infty} \frac{\{t\}}{t^\frac{4}{3}}dt - \int_{x^{1/4}}^{\infty} \frac{\{t\}}{t^\frac{4}{3}}dt\right] + \mathcal{O} \left(x^{1/4}\right) .\end{align*}
The improper integral $\int_{1}^{\infty} \frac{\{t\}}{t^\frac{4}{3}}dt$ exists since it is dominated by $\int_{1}^{\infty} \frac{1}{t^\frac{4}{3}}dt$. However, we know that \cite[p.232]{tenenbaum}
\begin{equation}\label{eqzeta1}
    \zeta\left(\frac{1}{3}\right) = -\frac{1}{2}- \frac{1}{3}\int_{1}^{\infty} \frac{\{t\}}{t^\frac{4}{3}}dt.
\end{equation}
Furthermore,
\begin{equation}\label{eqzeta2}
   0\leq \int_{x^{1/4}}^{\infty} \frac{\{t\}}{t^\frac{4}{3}}dt \leq \int_{x^{1/4}}^{\infty} \frac{1}{t^\frac{4}{3}}dt = \frac{3}{x^{\frac{1}{12}}},
\end{equation}
and hence
\begin{equation}\label{eq3.12}
    \int_{1}^{x^{1/4}} \{t\}\frac{f_1(t)}{t}dt= -6^{4/3} x^{1/3} \left(\zeta\left(\frac{1}{3}\right)+ \frac{1}{2}\right) + \mathcal{O}(x^{\frac{1}{4}})  = -6^{4/3}  x^{1/3} \left(\zeta\left(\frac{1}{3}\right)+ \frac{1}{2}\right) + \mathcal{O}(x^{\frac{1}{4}}) .
\end{equation}
Now, to compute $\int_{1}^{x^{1/4}} t\{t\}(g_1'(t)+ g_2'(t))dt$, we observe that
\begin{align*}
	g_1'(t)&= \frac{1}{3}\left(\frac{3x}{t^4}+\sqrt{\left(\frac{3x}{t^4}\right)^2-\frac{1}{27}}\right)^{-2/3} \frac{(-4)  (3x)}{t^5} \left[1+ \frac{\frac{3x}{t^4}} {\sqrt{\left(\frac{3x}{t^4}\right)^2-\frac{1}{27}}} \right] \\
	&= \frac{(-4x)}{t^5} \left(\frac{3x}{t^4}+\sqrt{\left(\frac{3x}{t^4}\right)^2-\frac{1}{27}}\right)^{1/3}  \left[ \frac{1} {\sqrt{\left(\frac{3x}{t^4}\right)^2-\frac{1}{27}}} \right],
\end{align*}
and 
\begin{align*}
	g_2'(t)&= \frac{1}{3}\left(\frac{3x}{t^4}- \sqrt{\left(\frac{3x}{t^4}\right)^2-\frac{1}{27}}\right)^{-2/3} \frac{(-4)  (3x)}{t^5} \left[1- \frac{\frac{3x}{t^4}} {\sqrt{\left(\frac{3x}{t^4}\right)^2-\frac{1}{27}}} \right] \\
	&= \frac{4x}{t^5} \left(\frac{3x}{t^4}- \sqrt{\left(\frac{3x}{t^4}\right)^2-\frac{1}{27}}\right)^{1/3}  \left[ \frac{1} {\sqrt{\left(\frac{3x}{t^4}\right)^2-\frac{1}{27}}} \right]. \\ 
\end{align*}
Thus, 
\begin{align}\label{eqA1}
\begin{split}
&g_1'(t)+ g_2'(t)= -\left(\frac{4x}{t^5} \right) \left[ \frac{1} {\sqrt{\left(\frac{3x}{t^4}\right)^2-\frac{1}{27}}} \right] \\
 & \times\left[\left(\frac{3x}{t^4}+\sqrt{\left(\frac{3x}{t^4}\right)^2-\frac{1}{27}}\right)^{1/3} - \left(\frac{3x}{t^4}-\sqrt{\left(\frac{3x}{t^4}\right)^2-\frac{1}{27}}\right)^{1/3}\right] \\ 
   &=  -\left(\frac{4x}{t^5} \right) \left[ \frac{ \left(\frac{3x}{t^4}\right)^{1/3} } {\sqrt{\left(\frac{3x}{t^4}\right)^2-\frac{1}{27}}} \right] \\
   &\times\left[\left(1+ \left(\frac{3x}{t^4}\right)^{-1} \sqrt{\left(\frac{3x}{t^4}\right)^2-\frac{1}{27}}\right)^{1/3} - \left( 1- \left(\frac{3x}{t^4}\right)^{-1} \sqrt{\left(\frac{3x}{t^4}\right)^2-\frac{1}{27}}\right)^{1/3}\right]
   \end{split}
  \end{align} 
Setting $\sqrt{u}:=  \left(\frac{3x}{t^4}\right)^{-1}\sqrt{\left(\frac{3x}{t^4}\right)^2-\frac{1}{27}}$, we get
\begin{align}\label{eqA2}
	 [(1+ \sqrt{u})^{1/3}-(1-\sqrt{u})^{1/3}]
	=  2 \sum_{m=0}^{\infty} \binom{1/3}{2m+1} \sqrt{u}^{2m+1}
	= 2 \sqrt{u} \left[\frac{1}{3}+ \frac{5}{81} u + \frac{22}{729}u^2+\cdots\right] 
\end{align}
Putting \eqref{eqA2} in \eqref{eqA1}, we infer that
\begin{align*}
	& g_1'(t)+ g_2'(t)= - 2^{10/3} 3^{-2/3} x^{1/3} t^{-7/3} + \mathcal{O} ( x^{-5/3} t^{17/3}).
\end{align*}
Therefore,
\begin{align}
\begin{split}\label{eq3.13}
    &\int_{1}^{x^{1/4}} t\{t\}(g_1'(t)+ g_2'(t))dt\\
    & = \int_{1}^{x^{1/4}} t\{t\}\left[-\frac{2^{10/3}\cdot 3^{-2/3}\cdot x^{1/3}}{t^{7/3}}\right]dt + \mathcal{O} \left( \int_{1}^{x^{1/4}} t\{t\} \left(x^{-5/3} t^{17/3}\right) dt \right)\\
    & = -2^{10/3}\cdot 3^{-2/3}\cdot x^{1/3} \int_{1}^{x^{1/4}}\frac{\{t\}}{t^{\frac{4}{3}}} dt + \mathcal{O} \left( x^{-5/3} \int_{1}^{x^{1/4}}  t^{20/3} dt \right)  \\
    & = 2^{10/3} \cdot 3^{1/3}\cdot x^{1/3} \left(\zeta\left(\frac{1}{3}\right)+ \frac{1}{2}\right) + \mathcal{O}(x^{\frac{1}{4}}).
    \end{split}
\end{align}
Substituting \eqref{eq3.12} and \eqref{eq3.13} in \eqref{eq3.9}, we get
\begin{equation}\label{eq3.14}
\int_{1}^{x^{1/4}} \{t\}f_1'(t)dt = 10 \cdot 2^{1/3}(3x)^{1/3}\left(\zeta\left(\frac{1}{3}\right)+ \frac{1}{2}\right) + \mathcal{O}(x^{\frac{1}{4}}).
\end{equation}
From \eqref{align3.6}, \eqref{eq3.6}, \eqref{eq3.8} and \eqref{eq3.14}, we infer 
\begin{equation}\label{ea3.16}
    \sum_{1\leq n\leq x^{\frac{1}{4}}}\lfloor{T_1\rfloor}= \left( \sqrt{3}G_1(0)-\frac{1}{2}\right) x^{\frac{1}{2}} + \left(10 \cdot 6^{\frac{1}{3}}\left(\zeta\left(\frac{1}{3} \right)+\frac{1}{2}\right)-\left(\frac{9}{2}\right)^{\frac{2}{3}}+2^{4/3}3^{1/3}\right)x^{\frac{1}{3}} + \mathcal{O}(x^{1/4}).
\end{equation}
Next consider
\begin{align}\label{align3.17}
\begin{split}
   \sum_{1\leq m\leq x^{\frac{1}{4}}}\lfloor{T_2\rfloor} &= \sum_{1\leq m\leq x^{\frac{1}{4}}}m\left(\left(\frac{3x}{2m^4}+\frac{1}{2}\sqrt{\left(\frac{3x}{m^4}\right)^2-\frac{1}{432}}\right)^{1/3} + \left(\frac{3x}{2m^4}-\frac{1}{2}\sqrt{\left(\frac{3x}{m^4}\right)^2-\frac{1}{432}}\right)^{1/3}\right)\\&-\sum_{1\leq m \leq x^{\frac{1}{4}}}\frac{m}{2}+ \mathcal{O}(x^{\frac{1}{4}})=\sum_{1\leq m\leq x^{\frac{1}{4}}} f_2(m)-{\frac{\lfloor x^{\frac{1}{4}} \rfloor (\lfloor x^{\frac{1}{4}} \rfloor+1)}{4}}+\mathcal{O}(x^{\frac{1}{4}})\\
   &= \sum_{1\leq m\leq x^{\frac{1}{4}}} f_2(m) -  \frac{1}{ 4 }  x^{\frac{1}{2}}+ \mathcal{O}(x^{\frac{1}{4}}), 
   \end{split}
\end{align} where
\begin{align*}
	f_2(m) &= m \left[\left(\frac{3x}{2m^4}+\frac{1}{2}\sqrt{\left(\frac{3x}{m^4}\right)^2-\frac{1}{432}}\right)^{1/3} + \left(\frac{3x}{2m^4}-\frac{1}{2}\sqrt{\left(\frac{3x}{m^4}\right)^2-\frac{1}{432}}\right)^{1/3}\right].
\end{align*}
By the Abel partial summation formula with $a_m=1, A(t) = \sum_{1\leq m\leq t}a_m$ and 
\begin{align*}
f_2(t) &= t \cdot 2^{-1/3} \left[\left(\frac{3x}{t^4}+ \sqrt{\left(\frac{3x}{t^4}\right)^2-\frac{1}{432}}\right)^{1/3} + \left(\frac{3x}{t^4}- \sqrt{\left(\frac{3x}{t^4}\right)^2-\frac{1}{432}}\right)^{1/3}\right],
\end{align*}
we get
\begin{align}\label{eq3.18b}
    \begin{split}
 \sum_{1\leq m\leq x^{\frac{1}{4}}} f_2(m)&= f_2(1)+\int_{1}^{x^{1/4}}f_2(t)dt+ \int_{1}^{x^{1/4}} \{t\}f_2'(t)dt+ \mathcal{O}(x^{1/4}). 
    \end{split}
\end{align}
Setting $\sqrt{r}:=  \left(\frac{3x}{t^4}\right)^{-1}\sqrt{\left(\frac{3x}{t^4}\right)^2-\frac{1}{432}}$, we deduce
\begin{align}\label{eqnewf2}
\begin{split}
	\frac{f_2(t)}{t} &=2^{-1/3} \left(\frac{3x}{t^4}\right)^{1/3}[(1+ \sqrt{r})^{1/3}+(1-\sqrt{r})^{1/3}]\\
	&= 2^{2/3} \left(\frac{3x}{t^4}\right)^{1/3}\left[1-\frac{1}{9}r-\frac{10}{243}r^2-\cdots\right] =  2 \left(\frac{3x}{t^4}\right)^{1/3} + \mathcal{O} \left( \left(3x\right)^{-5/3} t^8\right).
	\end{split}
\end{align}
From \eqref{eqnewf2} and \eqref{eq3.18b}, we get 
\begin{equation}\label{eq3.18}
 \sum_{1\leq m\leq x^{\frac{1}{4}}} f_2(m)= 2(3x)^{1/3}+\int_{1}^{x^{1/4}}f_2(t)dt+ \int_{1}^{x^{1/4}} \{t\}f_2'(t)dt+ \mathcal{O}(x^{1/4}). 
\end{equation}
Upon changing the variables $t\to \left(3x\right)^{\frac{1}{4}}t$, the first integral $\int_{1}^{x^{1/4}}f_2(t)dt$ is given by
\begin{align}\label{eq3.19}
\begin{split}
    &\int_{1}^{x^{\frac{1}{4}}}f_2(t)dt = 2^{-1/3} \left(3x\right)^{\frac{1}{2}}\\
    &\times \int_{\frac{1}{(3x)^{1/4}}}^{\left(\frac{1}{3}\right)^{1/4}}t\left[\left(\frac{1}{t^4}+\sqrt{\left(\frac{1}{t^4}\right)^2-\frac{1}{432}}\right)^{1/3} + \left(\frac{1}{t^4}-\sqrt{\left(\frac{1}{t^4}\right)^2-\frac{1}{432}}\right)^{1/3}\right]dt\\
    &=: 2^{-1/3} \left(3x\right)^{\frac{1}{2}}\int_{\frac{1}{(3x)^{1/4}}}^{\left(\frac{1}{3}\right)^{1/4}}G_2\left(\frac{1}{\left(3x\right)^{\frac{1}{4}}}\right)
\end{split}
\end{align}
where \[G_2(y):= \int_{y}^{\left(\frac{1}{3}\right)^{1/4}}t\left[\left(\frac{1}{t^4}+\sqrt{\left(\frac{1}{t^4}\right)^2-\frac{1}{432}}\right)^{1/3} + \left(\frac{1}{t^4}-\sqrt{\left(\frac{1}{t^4}\right)^2-\frac{1}{432}}\right)^{1/3}\right]dt.\]
Rewriting
\[G_2(y)= \int_{y}^{\left(\frac{1}{3}\right)^{1/4}} t^{-1/3} \left[\left(1+\sqrt{1-\frac{t^8}{432}}\right)^{1/3} + \left(1-\sqrt{1-\frac{t^8}{432}}\right)^{1/3}\right]dt\]
and then simplifying, 
\begin{align*}
G_2(0)& = \int_{0}^{y} t^{-1/3} \left[\left(1+\sqrt{1-\frac{t^8}{432}}\right)^{1/3} + \left(1-\sqrt{1-\frac{t^8}{432}}\right)^{1/3}\right]dt + G_2(y)\\
&= G_2(y)+ \int_{0}^{y}t^{-1/3} H\left(1-\frac{t^8}{432}\right)dt
\end{align*}
where $H(z)$ is given in \eqref{eqnewhz}. By using the binomial expansion \eqref{binom1}, we get
\begin{align*}
	G_2(0)= G_2(y)+ \int_{0}^{y} \left( 2^{1/3}t^{-1/3} + \mathcal{O} \left(t^{23/3}\right) \right) dt
\end{align*}
 and hence
\[G_2(y) = G_2(0) - 3\cdot 2^{-\frac{2}{3}} y^{\frac{2}{3}} +\mathcal{O}(y^{\frac{26}{3}}).\]
Thus,
\begin{align}
\begin{split}\label{eq3.20}
    \int_{1}^{x^{\frac{1}{4}}}f_2(t)dt &= 2^{-1/3} \left(3x\right)^{\frac{1}{2}} \left[G_2(0) - 3\cdot 2^{-\frac{2}{3}} \left(\frac{1}{\left(3x\right)^{\frac{1}{4}}}\right)^{\frac{2}{3}} +\mathcal{O} \left(\left(\frac{1}{\left(3x\right)^{\frac{1}{4}}}\right)^{\frac{26}{3}}\right)\right]\\
    &=  2^{-1/3} \left(3x\right)^{\frac{1}{2}}G_2(0) - 3^{\frac{4}{3}}\cdot 2^{-1} x^{\frac{1}{3}} + \mathcal{O}(x^{-\frac{5}{3}}).
    \end{split}
\end{align}
To compute $\int_{1}^{x^{1/4}} \{t\}f_2'(t)dt$, set $f_2(t) = t(g_3(t)+g_4(t))$ where
\begin{align*}
	g_3(t)= 2^{-1/3} \left(\frac{3x}{t^4}+ \sqrt{\left(\frac{3x}{t^4}\right)^2-\frac{1}{432}}\right)^{1/3}, \quad g_4(t)= 2^{-1/3} \left(\frac{3x}{t^4}- \sqrt{\left(\frac{3x}{t^4}\right)^2-\frac{1}{432}}\right)^{1/3}.
\end{align*}
Then
\begin{equation}\label{eq3.21B}
	\int_{1}^{x^{1/4}} \{t\}f_2'(t)dt = \int_{1}^{x^{1/4}} \{t\}\frac{f_2(t)}{t}dt+ \int_{1}^{x^{1/4}} t\{t\}(g_3'(t)+ g_4'(t))dt.
\end{equation}

Hence,
\begin{align*}
	\int_{1}^{x^{1/4}} \{t\}\frac{f_2(t)}{t}dt& = 2(3x)^{1/3}\int_{1}^{x^{1/4}} \frac{\{t\}}{t^\frac{4}{3}}dt + \mathcal{O} \left( \left(3x\right)^{-5/3} \int_{1}^{x^{1/4}} \{t\} t^{20/3} dt \right) \\
	& = 2(3x)^{1/3}\left[\int_{1}^{\infty} \frac{\{t\}}{t^\frac{4}{3}}dt - \int_{x^{1/4}}^{\infty} \frac{\{t\}}{t^\frac{4}{3}}dt\right] + \mathcal{O} \left(x^{1/4}\right) .\end{align*}
Using \eqref{eqzeta1} and \eqref{eqzeta2}, we infer
\begin{align}\label{eq3.12B}
\begin{split}
	\int_{1}^{x^{1/4}} \{t\}\frac{f_2(t)}{t}dt&= (-6)(3x)^{1/3} \left(\zeta\left(\frac{1}{3}\right)+ \frac{1}{2}\right) + \mathcal{O}(x^{\frac{1}{4}}).
	\end{split}
\end{align}
Note that
\begin{align*}
	g_3'(t)&= \frac{2^{-1/3}}{3}\left(\frac{3x}{t^4}+\sqrt{\left(\frac{3x}{t^4}\right)^2-\frac{1}{432}}\right)^{-2/3} \frac{(-4)  (3x)}{t^5} \left[1+ \frac{\frac{3x}{t^4}} {\sqrt{\left(\frac{3x}{t^4}\right)^2-\frac{1}{432}}} \right] \\
	&= \frac{(-2^{5/3}x)}{t^5} \left(\frac{3x}{t^4}+\sqrt{\left(\frac{3x}{t^4}\right)^2-\frac{1}{432}}\right)^{1/3}  \left[ \frac{1} {\sqrt{\left(\frac{3x}{t^4}\right)^2-\frac{1}{432}}} \right],
	\end{align*}
and
\begin{align*}
	g_4'(t)&= \frac{2^{-1/3}}{3}\left(\frac{3x}{t^4}- \sqrt{\left(\frac{3x}{t^4}\right)^2-\frac{1}{432}}\right)^{-2/3} \frac{(-4)  (3x)}{t^5} \left[1- \frac{\frac{3x}{t^4}} {\sqrt{\left(\frac{3x}{t^4}\right)^2-\frac{1}{432}}} \right] \\
	&= \frac{2^{5/3}x}{t^5} \left(\frac{3x}{t^4}- \sqrt{\left(\frac{3x}{t^4}\right)^2-\frac{1}{432}}\right)^{1/3}  \left[ \frac{1} {\sqrt{\left(\frac{3x}{t^4}\right)^2-\frac{1}{432}}} \right]. \\ 
\end{align*}
Now, we have
\begin{align}\label{eqB1}
\begin{split}
	& g_3'(t)+ g_4'(t)=  -\left(\frac{2^{5/3}x}{t^5} \right) \left[ \frac{ \left(\frac{3x}{t^4}\right)^{1/3} } {\sqrt{\left(\frac{3x}{t^4}\right)^2-\frac{1}{432}}} \right] \\
	&\times \left[\left(1+ \left(\frac{3x}{t^4}\right)^{-1} \sqrt{\left(\frac{3x}{t^4}\right)^2-\frac{1}{432}}\right)^{1/3} - \left( 1- \left(\frac{3x}{t^4}\right)^{-1} \sqrt{\left(\frac{3x}{t^4}\right)^2-\frac{1}{432}}\right)^{1/3}\right].
	\end{split}
\end{align} 
Setting $\sqrt{v}:=  \left(\frac{3x}{t^4}\right)^{-1}\sqrt{\left(\frac{3x}{t^4}\right)^2-\frac{1}{432}}$ and using
\eqref{eqA2} in \eqref{eqB1}, we get
\begin{align*}
	& g_3'(t)+ g_4'(t)= - 2^3 3^{-2/3} x^{1/3} t^{-7/3} + \mathcal{O} ( x^{-5/3} t^{17/3})
\end{align*}
Therefore,
\begin{align}
	\begin{split}\label{eq3.13B}
		&\int_{1}^{x^{1/4}} t\{t\}(g_3'(t)+ g_4'(t))dt= \int_{1}^{x^{1/4}} t\{t\}\left[-\frac{2^{3}\cdot 3^{-2/3}\cdot x^{1/3}}{t^{7/3}}\right]dt \\
		&+ \mathcal{O} \left( \int_{1}^{x^{1/4}} t\{t\} \left( x^{-5/3} t^{17/3}\right) dt \right)\\
		& = -2^{3}\cdot 3^{-2/3}\cdot x^{1/3} \int_{1}^{x^{1/4}}\frac{\{t\}}{t^{\frac{4}{3}}} dt + \mathcal{O} \left( x^{-5/3} \int_{1}^{x^{1/4}}  t^{20/3} dt \right)  \\
		& = 2^{3} \cdot 3^{1/3}\cdot x^{1/3} \left(\zeta\left(\frac{1}{3}\right)+ \frac{1}{2}\right) + \mathcal{O}(x^{\frac{1}{4}}).
	\end{split}
\end{align}
Substituting \eqref{eq3.12B} and \eqref{eq3.13B} in \eqref{eq3.21B}, we get
\begin{equation}\label{eq3.14B}
	 \int_{1}^{x^{1/4}} \{t\}f_2'(t)dt =  2(3x)^{1/3}\left(\zeta\left(\frac{1}{3}\right)+ \frac{1}{2}\right) + \mathcal{O}(x^{\frac{1}{4}}).
\end{equation}
From \eqref{align3.17}, \eqref{eq3.18}, \eqref{eq3.20} and \eqref{eq3.14B}, we infer 
\begin{equation}\label{eq3.25}
    \sum_{1\leq m\leq x^{\frac{1}{4}}}\lfloor{T_2\rfloor}= \left( 2^{-1/3}3^{1/2} G_2(0)-\frac{1}{4}\right)x^{\frac{1}{2}} + \left(2\cdot 3^{\frac{1}{3}}\left(\zeta\left(\frac{1}{3} \right)+\frac{1}{2}\right)-3^{\frac{4}{3}}2^{-1}+2\cdot3^{1/3}\right)x^{\frac{1}{3}} + \mathcal{O}(x^{1/4}).
\end{equation}
Finally, from \eqref{align3.2} we have
\begin{align*}
\sum_{1\leq n\leq x}\varrho_{\mathfrak{so}(5)}(n) &= \left(\sqrt{3}G_1(0)+2^{-1/3}\sqrt{3}G_2(0)-\frac{3}{4}\right)x^{\frac{1}{2}} \\
& + \left( 3^{1/3} \left(\zeta\left(\frac{1}{3} \right)+\frac{1}{2}\right) (10 \cdot 2^{1/3}+2) - 3^{4/3} (2^{-2/3}+2^{-1})+3^{\frac{1}{3}}(2+2^{\frac{4}{3}} ) \right)x^{\frac{1}{3}} + \mathcal{O}(x^{\frac{1}{4}}).
\end{align*}
By comparing with \eqref{eq3.1}, we conclude the proof of Theorem \ref{thm1}.
\begin{rmk}
We deduce from the proof of Theorem \ref{thm1} that
\begin{align*}
&\int_{0}^{\left(\frac{1}{3}\right)^{1/4}}t\left[\left(\frac{1}{t^4}+\sqrt{\left(\frac{1}{t^4}\right)^2-\frac{1}{27}}\right)^{1/3} + \left(\frac{1}{t^4}-\sqrt{\left(\frac{1}{t^4}\right)^2-\frac{1}{27}}\right)^{1/3}\right]dt\\
&+\int_{0}^{\left(\frac{1}{3}\right)^{1/4}}t\left[\left(\frac{1}{2t^4}+\sqrt{\left(\frac{1}{2t^4}\right)^2-\frac{1}{1728}}\right)^{1/3} + \left(\frac{1}{2t^4}-\sqrt{\left(\frac{1}{2t^4}\right)^2-\frac{1}{1728}}\right)^{1/3}\right]dt\\
&= \frac{\sqrt{3}}{4}+ \frac{\sqrt{3}\Gamma(\frac{1}{4})^2}{4\sqrt{\pi}}.
\end{align*}
Since the right hand side of the above equality is a transcendental number by Nesterenko's result, we can deduce that at least one of the integral in left hand side is a transcendental number.
\end{rmk}
\begin{rmk}
It would be interesting to consider deeper techniques to improve the error term in Theorem \ref{thm1}.
\end{rmk}

\section*{Acknowledgements}
The first author has received funding support from the Anusandhan National Research Foundation under the Start-up Research Grant (File No. SRG/2023/002255). S.S.R. work is supported by a grant from the ANRF (File No.:CRG/2022/000268).


\begin{thebibliography}{9999}



\bibitem{Bringmann2023}
K. Bringmann and J. Franke1, An asymptotic formula for the number of $n$-dimensional representations of $\mathfrak{su}(3)$, \textit{Rev. Mat. Iberoam.} \textbf{39}(5) (2023), 1599-1640.

\bibitem{Bridges2024}
W. Bridges, B. Brindle, K. Bringmann and J. Franke1, Asymptotic expansions for partitions generated by infinite products, \textit{Math. Ann.} \textbf{390} (2024), 2593-2632.

\bibitem{Bridges-Bringmann2024}
W. Bridges, K. Bringmann and J. Franke1, On the number of irreducible representaions of $\mathfrak{su}(3)$, \textit{Acta Arith.} \textbf{215} (2024), 65-71.

\bibitem{gunnels-sczech}
P. E. Gunnels, R. Sczech, Evaluation of Dedekind sums, Eisenstein cocycles, and special values of $L$-functions. \textit{Duke. Math. J.} \textbf{118} (2003), 229--260.

\bibitem{hardy-ramanujan}
G. H. Hardy, S. Ramanujan.
Asymptotic formulae in combinatory analysis.
\textit{Proc. London Math. Soc.} (2) 17 (1918), 75--115.

\bibitem{hall}
B. C. Hall. Lie Groups, Lie Algebras, and Representations: An Elementary Introduction. Springer, 2004.


\bibitem{matsumoto}
K. Matsumoto. On the analytic continuation of various multiple zeta-functions. In: Number Theory for the Millenium II: Proceedings of the Millennial Conference on Number Theory (Urbana-Champaign, USA, 2000), pp.~417--440. Eds. B. C. Berndt, N. Boston, H. G. Diamond, A. J. Hildebrand. AK~Peters, 2002.

\bibitem{matsumoto-bonner}
K. Matsumoto. On Mordell-Tornheim and other multiple zeta functions. In: 
Proceedings of the Session in Analytic Number Theory and Diophantine Equations, Eds. D. R. Heath-Brown, B. Z. Moroz.
Bonner Mathematische Schriften vol.~360, Univ. Bonn, 2003.

\bibitem{matsumoto-semisimple1}
K. Matsumoto, H. Tsumura.
On Witten multiple zeta-functions associated with semisimple Lie algebras I.
\textit{Ann. Inst. Fourier} \textbf{56} (2006), 1457--1504.

\bibitem{knapp}
A. W. Knapp, Representation theory of semisimple groups, Princeton University Press, Princeton and Oxford (1986).


\bibitem{matsumoto-semisimple2}
Y. Komori, K. Matsumoto, H. Tsumura.
On Witten multiple zeta-functions associated with semisimple Lie algebras II.
\textit{J. Math. Soc. Japan} \textbf{62} (2010), 355--394.

\bibitem{matsumoto-semisimple3}
Y. Komori, K. Matsumoto, H. Tsumura.
On Witten multiple zeta-functions associated with semisimple Lie algebras III.
In: Multiple Dirichlet Series, $L$-functions and Automorphic Forms, pp.~223--286. Eds. D. Bump, S. Friedberg, D. Goldfeld. Birkh\"auser, 2012.

\bibitem{matsumoto-semisimple4}
Y. Komori, K. Matsumoto, H. Tsumura.
On Witten multiple zeta-functions associated with semisimple Lie algebras IV.
\textit{Glasgow Math. J.} \textbf{53} (2011), 185--206.


\bibitem{mordell}
L. J. Mordell, On the evaluation of some multiple series. \textit{J. London Math. Soc.} \textbf{33}, 1958, 368--371.

\bibitem{paris-kaminski}
R. B. Paris, D. Kaminski. Asymptotics and Mellin-Barnes Integrals. Encyclopedia of Mathematics and its Applications, Vol. 85. Cambridge University Press, 2001.

\bibitem{murty}
M. Ram Murty, Problems in Analytic Number Theory, Springer-Verlag New York (2008).



\bibitem{romik2017}
D. Romik, On the number of $n$-dimensional representations of $SU(3)$, the Bernoulli numbers, and the Witten zeta function, \textit{Acta Arith.} \textbf{180} (2017), 111-159.

\bibitem{subbarao-etal}
M. V. Subbarao, R. Sitaramachandrarao, On some infinite series of L. J. Mordell and their analogues. 
\textit{Pacific J. Math.} \textbf{19} (1985), 245--255.

\bibitem{tenenbaum}G. Tenenbaum, Introduction to Analytic and Probabilistic Number Theory, Graduate Studies in Mathematics, AMS {\bf 163} (2008).

\bibitem{tornheim}
L. Tornheim, Harmonic double series. \textit{Amer. J. Math.} \textbf{72} (1950), 303--314.


\bibitem{witten}
E. Witten, On quantum gauge theories in two dimensions. \textit{Commun. Math. Phys.} \textbf{141} (1991), 153--209.

\bibitem{zagier}
D. Zagier, Values of zeta functions and their applications. In: First European Congress of Mathematics, Volume II, Prog. Math., 120, Birkh\"auser, Basel, 1994, 497–512.

\end{thebibliography}
\end{document}